\documentclass[12pt,amssymb, amsfonts]{article}
\usepackage{graphicx, amssymb}
\usepackage{algorithmic}
\usepackage{algorithm}
\usepackage{authblk}

\newtheorem{theorem}{Theorem}[section]
\newtheorem{corollary}[theorem]{Corollary}

\newtheorem{lemma}[theorem]{Lemma}
\newtheorem{example}[theorem]{Example}
\newtheorem{remark}[theorem]{Remark}

\title{A Note on the Construction of  Complex and Quaternionic
       Vector Fields on Spheres
       \thanks{Published in the Journal of  Mathematical Notes,  93(1)  (2013),  104-110}}

\author{Mohammad Obiedat}
\affil{Department of Science, Technology, and Mathematics, Gallaudet University,
800 Florida Avenue NE,  Washington,  DC 20002-3695, USA\thanks{Email address: mohammad.obiedat@gallaudet.edu}}

\begin{document}
\maketitle

\begin{abstract}
A relationship between  real, complex, and quaternionic vector fields on spheres is given by using   a relationship between the corresponding standard inner products.  The number of linearly independent complex vector fields on the standard $(4n-1)$-sphere is shown to be twice the number of linearly independent quaternionic vector fields  plus $d$, where $ d = 1 \mbox{ or } 3$.

\vspace{0.5cm}
\end{abstract}

\section{Introduction}

Let ${\mathbb F}$  be one of  the associative division algebras
over $\mathbb{R}$, namely, the real numbers  ${\mathbb R}$, the
complex numbers ${\mathbb C} = \{ a+bi:a,b\in {\mathbb R} \}$, or
the quaternionic numbers ${\mathbb H} = \{a+bi+cj+dk:a,b,c,d\in
{\mathbb R}\}$, with the usual norm and conjugation operations.
For each $n \in  {\mathbb N}$, we consider ${\mathbb F}^n$ as a
right inner product space over ${\mathbb F}$ with  the usual
multiplication and the standard ${\mathbb F}$-inner product
$\langle|\rangle_{\mathbb F}:{\mathbb F}^n\times {\mathbb F}^n\rightarrow
{\mathbb F}$ given by
$\langle(x_1,\ldots,x_n)|(y_1,\ldots,y_n)\rangle_{\mathbb F}
=\sum_{m=1}^{n}\bar{x}_m y_m$. Let  $S({\mathbb F}^n)=\{x\in
{\mathbb F}^n:\langle x|x\rangle_{\mathbb F} = 1\}$ be the standard unit sphere in
${\mathbb F}^n$. An ${\mathbb F}$-vector field on $S({\mathbb
F}^n)$ is a continuous function $v:S({\mathbb F}^n)\rightarrow
{\mathbb F}^n- \{0\}$ such that $\langle x|v(x)\rangle_{\mathbb F}=0$ for each
$x\in S({\mathbb F}^n)$. Given $m$ such ${\mathbb F}$-vector  fields 
 $v_1,\ldots,v_m$   on $S({\mathbb F}^n)$,
we say that they are linearly independent if the vectors $v_1(x),\ldots,v_m(x)$ are linearly independent over ${\mathbb F}$ for all $x\in S({\mathbb F}^n)$.
Let $\rho^{\mathbb F}({\mathbb F}^n)$ denote  the maximal number of linearly
independent ${\mathbb F}$-vector fields on $S({\mathbb F}^n)$. The
vector fields on spheres problem has two sides: the first side is
what we call the maximal number problem, namely, the computation
of $\rho^{\mathbb F}({\mathbb F}^n)$, the second side is what we
call the construction problem, namely, the actual construction of
$\rho^{\mathbb F}({\mathbb F}^n)$ linearly independent
${\mathbb F}$-vector fields on $S({\mathbb F}^n)$. The roots of this
classical problem goes back to the hairy ball theorem, the
parallelizability of spheres, and  the classification of division
algebras over ${\mathbb R}$. For background materials on this
problem, we refer the reader to \cite{MPS80, Tho69}.

If $n$ is odd, then one can easily show that
$\rho^{\mathbb F}({\mathbb F}^n)=0$. Therefore, unless otherwise indicated,
we will assume that $n$ is an even positive integer.
An explicit formula for computing $\rho^{\mathbb R}({\mathbb R}^n)$ was
given by Adams \cite{Ada62} in 1962,  when he showed that
$\rho^{\mathbb R}({\mathbb R}^n) = \rho(n)-1$, where
$\rho(n)$ be the Radon-Hurwitz function.
While there are no explicit  formulas for computing
$\rho^{\mathbb C}({\mathbb C}^n)$ and $\rho^{\mathbb H}({\mathbb H}^n)$,
they still can be obtained by using the work of  Adams and Walker
\cite{AW65} in 1965 for  the complex case, and the work of
Sigrist and Suter \cite{SS73} in 1973 for the quaternionic case.
More direct  formulas for  computing
$\rho^{\mathbb C}({\mathbb C}^n)$ and $\rho^{\mathbb H}({\mathbb H}^n)$
are given in Section~3 of this paper.

Our main goal in this paper is to draw more attention to the
importance of the construction problem of vector fields on spheres.
There are several known  methods that explicitly give
$\rho^{\mathbb R}({\mathbb R}^n)$ linearly independent real vector fields on
$S({\mathbb R}^n)$ (e.g., see \cite{Ogn08, Zve68}).
The situation is completely different with  the construction of
complex and quaternionic vector fields;  there is no explicit construction
that  gives two or more linearly independent complex vector fields on
$S({\mathbb C}^n)$, and there is no known construction that
gives even a single quaternionic vector field on $S({\mathbb H}^n)$.
In addition to their self importance, the actual construction of complex
and quaternionic vector fields on spheres  might lead to the solution
of several open problems in the equivairant complex and quaternionic
vector fields on  spheres (see \cite{Obi06, Ond01}).

The layout of this paper is as follows.
In Section~2, we present a  relationship between real, complex, and
quaternionic standard  inner products. We show how such a relationship
leads to  a relationship between corresponding real, complex,
and quaternionic vector fields on spheres. More specifically,
we  show how $m$ linearly independent  quaternionic vector fields can
be used to obtain $2m$ linearly independent complex vector fields,
and how $m$ linearly independent  complex  vector fields can be used to obtain
$2m$ linearly independent real  vector fields.
In Section 3, we provide direct formulas for computing
$\rho^{\mathbb C}({\mathbb C}^n)$ and $\rho^{\mathbb H}({\mathbb H}^n)$, and show that
$\rho^{\mathbb C}({\mathbb C}^{2n}) = 2\rho^{\mathbb H}({\mathbb H}^n)+d$
where $d = 1\mbox{ or } 3$.
In Section~4, we give necessary and sufficient conditions on
linearly independent real (respectively,  complex ) vector
fields to be linearly independent complex  or quaternionic
(respectively,  quaternionic) vector fields.

\section{A relationship between real, complex, and quaternionic vector fields on spheres}
For any positive integer $n$, let
$r_{\mathbb C}:{\mathbb C}^n\rightarrow {\mathbb R}^{2n}$ be the realification function from ${\mathbb C}^n$ to ${\mathbb R}^{2n}$ defined by
$r_{\mathbb C}(a_1+b_1i,\ldots,a_n+b_ni)=(a_1,b_1,\ldots,a_n,b_n)$,
$c_{\mathbb H}:{\mathbb H}^n\rightarrow {\mathbb C}^{2n}$ be the complexification function from ${\mathbb H}^n$ to ${\mathbb C}^{2n}$ defined by
$c_{\mathbb H}(a_1+b_1i+c_1j+d_1k,\ldots,a_n+b_ni+c_nj+d_nk)=
(a_1+b_1i,d_1+c_1i,\ldots,a_n+b_ni,d_n+c_ni)$, and
$r_{\mathbb H}:{\mathbb H}^n\rightarrow {\mathbb R}^{4n}$ be the realification function from ${\mathbb H}^n$ to ${\mathbb R}^{4n}$ defined by
$r_{\mathbb H} = r_{\mathbb C}\circ c_{\mathbb H}$.
For $t\in {\mathbb F}$, let $\alpha_t:{\mathbb F}^n\rightarrow {\mathbb F}^n$
be the right multiplication by $t$ , i.e.,
$\alpha_t(x)=x\cdot t$ for each $x \in {\mathbb F}^n$.
The proof of the following lemma is straightforward.

\begin{lemma} \label{lem:Composition}
Let $s\in {\mathbb C}$ and $t\in {\mathbb H}$. Then
\begin{enumerate}
\item[\textup{(i)}] $c_{\mathbb H}\circ \alpha_s = \alpha_s \circ c_{\mathbb H}$.
\item[\textup{(ii)}] $r_{\mathbb C}\circ \alpha_s \circ c_{\mathbb H}=
r_{\mathbb H}\circ \alpha_s$.
\item[\textup{(iii)}] $r_{\mathbb C}\circ \alpha_s \circ c_{\mathbb H}\circ \alpha_t=
r_{\mathbb H}\circ \alpha_{ts}$.
\end{enumerate}
\end{lemma}

In the following theorem, we give a relationship between real, complex, and
quaternionic standard inner products on ${\mathbb F}^n$.

\begin{theorem} \label{theorem:InnerProduct}
Let $x,y\in {\mathbb C}^n$ and $v,w\in {\mathbb H}^n$. Then
\begin{enumerate}
\item[\textup{(i)}]$\langle x|y\rangle_{\mathbb C}=\langle r_{\mathbb C}(x)|r_{\mathbb C}(y)\rangle_{\mathbb R}
- \langle r_{\mathbb C}(x)|r_{\mathbb C}\circ \alpha_i(y)\rangle_{\mathbb R}i$.
\item[\textup{(ii)}] $\langle v|w\rangle_{\mathbb H}=
\langle c_{\mathbb H}(v)|c_{\mathbb H}(w)\rangle_{\mathbb C}- \langle c_{\mathbb H}(v)|c_{\mathbb H}\circ \alpha_j(w)\rangle_{\mathbb C}j$.
\item[\textup{(iii)}]$\langle v|w\rangle_{\mathbb H}=
\langle r_{\mathbb H}(v)|r_{\mathbb H}(w)\rangle_{\mathbb R}
- \langle r_{\mathbb H}(v)|r_{\mathbb H}\circ \alpha_i(w)\rangle_{\mathbb R}i\\
-\langle r_{\mathbb H}(v)|r_{\mathbb H}\circ \alpha_j(w)\rangle_{\mathbb R}j
-\langle r_{\mathbb H}(v)|r_{\mathbb H}\circ \alpha_k(w)\rangle_{\mathbb R}k$.
\end{enumerate}
\end{theorem}
\textbf{Proof.} Without loss of generality, we can assume that $n = 1$. Then  the result  follows by using  Lemma~{\ref{lem:Composition}} and the definition of
the standard inner products.\\

In the following theorem, we give a relationship between real, complex, and
quaternionic vector fields on spheres.

\begin{theorem} \label{theorem:RelationshipVectorFields}
Let $v:S({\mathbb C}^n)\rightarrow {\mathbb C}^n- \{0\}$ and
$w:S({\mathbb H}^n)\rightarrow {\mathbb H}^n- \{0\}$
be two continuous functions. Then
\begin{enumerate}
\item[\textup{(i)}] $v$ is a complex vector field on $S({\mathbb C}^n)$ if and only if
$r_{\mathbb C}\circ v\circ r_{{\mathbb C}}^{-1}$ and
$r_{\mathbb C}\circ \alpha_i \circ v\circ r_{{\mathbb C}}^{-1}$ are real vector fields on $S({\mathbb R}^{2n})$.
\item[\textup{(ii)}] $w$ is a quaterionic vector field on $S({\mathbb H}^n)$ if and only if
$c_{\mathbb H}\circ w\circ c_{{\mathbb H}}^{-1}$ and
$c_{\mathbb H}\circ \alpha_j \circ w\circ c_{{\mathbb H}}^{-1}$ are complex vector fields on $S({\mathbb C}^{2n})$.
\item[\textup{(iii)}]$w$ is a quaterionic vector field on $S({\mathbb H}^n)$ if and only if
$r_{\mathbb H}\circ w\circ r_{{\mathbb H}}^{-1}$ and
$r_{\mathbb H}\circ \alpha_t \circ w\circ r_{{\mathbb H}}^{-1}$, where
$t\in \{i,j,k\}$, are real  vector fields on $S({\mathbb R}^{4n})$.
\end{enumerate}
\end{theorem}
\textbf{Proof.} (i) Let $x\in S({\mathbb R}^{2n})$.
By Theorem~{\ref{theorem:InnerProduct} (i)},
$\langle r_{\mathbb C}^{-1}(x)|v(r_{\mathbb C}^{-1}(x))\rangle_{\mathbb C} = 0$ if and only if
$\langle x|r_{\mathbb C}\circ v\circ r_{{\mathbb C}}^{-1}(x)\rangle_{\mathbb R}= 0$ and
$\langle x|r_{\mathbb C}\circ \alpha_i\circ v\circ r_{{\mathbb C}}^{-1}(x)\rangle_{\mathbb R} = 0$.
The result follows. Similarly, one can prove (ii) and (iii) by
using Theorem~{\ref{theorem:InnerProduct} (ii) and (iii)}.

\begin{example} \textup{\label{exa:ComplexVectorField}
For each $n\geqslant 2$, $v:S({\mathbb C}^n)\rightarrow {\mathbb C}^n- \{0\}$ such
that $v(x_1+ix_2,\ldots,x_{2n-1}+ix_{2n})
= (-x_3+ix_4,x_1-ix_2,\ldots,-x_{2n-1}+ix_{2n},x_{2n-3}-ix_{2n-2})$ is a complex vector field on $S({\mathbb C}^n)$. By Theorem~{\ref{theorem:RelationshipVectorFields} (i)},
$r_{\mathbb C}\circ v\circ r_{{\mathbb C}}^{-1}$ and
$r_{\mathbb C}\circ \alpha_i \circ v\circ r_{{\mathbb C}}^{-1}$ are real vector fields on
$S({\mathbb R}^{2n})$. $r_{\mathbb C}\circ v\circ r_{{\mathbb C}}^{-1}$ is given by
$(x_1, \ldots, x_{2n})\mapsto
(-x_3, x_4, x_1, -x_2, \ldots,-x_{2n-1}, x_{2n}, x_{2n-3}, -x_{2n-2})$
and $r_{\mathbb C}\circ \alpha_i \circ v\circ r_{{\mathbb C}}^{-1}$  is given by
$(x_1,\ldots, x_{2n})\mapsto
(-x_4, -x_3, x_2, x_1, \ldots,-x_{2n}, -x_{2n-1}, x_{2n-2}, x_{2n-3})$.}
\end{example}

\begin{theorem} \label{theorem:LinearlyIndependentVectorFields}
Let $v_1,\ldots,v_m:S({\mathbb C}^n)\rightarrow {\mathbb C}^n- \{0\}$ and
$w_1,\ldots,w_m:S({\mathbb H}^n)\rightarrow {\mathbb H}^n- \{0\}$
be continuous functions. Then
\begin{enumerate}
\item[\textup{(i)}] $v_1,\ldots,v_m$ are linearly independent  complex vector field on $S({\mathbb C}^n)$ if and only if
$r_{\mathbb C}\circ v_1\circ r_{{\mathbb C}}^{-1},
r_{\mathbb C}\circ \alpha_i \circ v_1\circ r_{{\mathbb C}}^{-1},
\ldots, r_{\mathbb C}\circ v_m\circ r_{{\mathbb C}}^{-1},
r_{\mathbb C}\circ \alpha_i \circ v_m\circ r_{{\mathbb C}}^{-1}$ are linearly independent real vector fields on $S({\mathbb R}^{2n})$.
\item[\textup{(ii)}]  $w_1,\ldots,w_m$ are linearly independent  quaternionic vector
field on $S({\mathbb H}^n)$ if and only if
$c_{\mathbb H}\circ w_1\circ c_{{\mathbb H}}^{-1},
c_{\mathbb H}\circ \alpha_j \circ w_1\circ r_{{\mathbb H}}^{-1},
\ldots, c_{\mathbb H}\circ w_m\circ c_{{\mathbb H}}^{-1},
c_{\mathbb H}\circ \alpha_j \circ w_m\circ c_{{\mathbb H}}^{-1}$ are linearly independent complex vector fields on $S({\mathbb C}^{2n})$.
\item[\textup{(iii)}]$w_1,\ldots,w_m$ are linearly independent  quaternionic vector
field on  $S({\mathbb H}^n)$ if and only if
$r_{\mathbb H}\circ w_1\circ r_{{\mathbb H}}^{-1},
r_{\mathbb H}\circ \alpha_t \circ w_1\circ r_{{\mathbb H}}^{-1},
\ldots, r_{\mathbb H}\circ w_m\circ r_{{\mathbb H}}^{-1},
r_{\mathbb H}\circ \alpha_t \circ w_m\circ r_{{\mathbb H}}^{-1}$, where $t\in\{i,j,k\}$, are linearly independent real vector fields on $S({\mathbb R}^{4n})$.
\end{enumerate}
\end{theorem}
\textbf{Proof.} (i) By Theorem~{\ref{theorem:RelationshipVectorFields} (i)}, $v_1,\ldots,v_m$ are  complex vector fields on $S({\mathbb C}^n)$
if and only if $r_{\mathbb C}\circ v_1\circ r_{{\mathbb C}}^{-1},
r_{\mathbb C}\circ \alpha_i \circ v_1\circ r_{{\mathbb C}}^{-1},
\ldots, r_{\mathbb C}\circ v_m\circ r_{{\mathbb C}}^{-1},
r_{\mathbb C}\circ \alpha_i \circ v_m\circ r_{{\mathbb C}}^{-1}$
are real vector fields on $S({\mathbb R}^{2n})$. Let $x\in S({\mathbb R}^{2n}), y=r_{{\mathbb C}}^{-1}(x)$, and $a_1,b_1,\ldots, a_m,b_m\in {\mathbb R}$.
Then  $(r_{\mathbb C}\circ v_1\circ r_{{\mathbb C}}^{-1}(x))a_1 +
(r_{\mathbb C}\circ \alpha_i \circ v_1\circ r_{{\mathbb C}}^{-1}(x))b_1 +
\cdots + (r_{\mathbb C}\circ v_m\circ r_{{\mathbb C}}^{-1}(x))a_m +
(r_{\mathbb C}\circ \alpha_i \circ v_m\circ r_{{\mathbb C}}^{-1}(x))b_m = 0$  if and only if $v_1(y)(a_1+b_1i)+\cdots+v_m(y)(a_m+b_mi) = 0$. The result follows.
(ii) and (iii) are similar to (i).

\begin{corollary} \label{corollary:RelationshipNumberOfVectorFields}
  $\rho^{\mathbb R}({\mathbb R}^{2n})\geqslant2\rho^{\mathbb C}({\mathbb C}^n)$
  and
  $\rho^{\mathbb C}({\mathbb C}^{2n})\geqslant 2\rho^{\mathbb H}({\mathbb H}^n)$.
\end{corollary}
\begin{remark} \label{remark:OrthonormalVectorFields}
\textup{A unit ${\mathbb F}$-vector field on $S({\mathbb F}^n)$ is an ${\mathbb F}$-vector field on $S({\mathbb F}^n)$ whose image is a subset of $S({\mathbb F}^n)$. $m$ such unit vector fields $v_1,\ldots,v_m$ are orthonormal if
$\langle v_s(x)|v_t(x)\rangle_{\mathbb F}=\delta_{st}$ for each $s,t\in\{1,\ldots,m\}$ and all $x\in S({\mathbb F}^n)$. Observe that since the Gram-Schmidt orthonormalization process is continuous, then one can convert $m$ linearly independent
${\mathbb F}$-vector fields on $S({\mathbb F}^n)$ to $m$ orthonormal ${\mathbb F}$-vector fields on $S({\mathbb F}^n)$. Consequently, $\rho^{\mathbb F}({\mathbb F}^n)$ is equal to the maximal number of orthonormal ${\mathbb F}$-vector fields on $S({\mathbb F}^n)$.
Also, one can use Theorem~\ref{theorem:InnerProduct} and basic properties of standard
inner products to prove that Theorem~\ref{theorem:RelationshipVectorFields} remains valid if we replace ``vector field'' by ``unit vector field'' , and Theorem~\ref{theorem:LinearlyIndependentVectorFields} remains valid if
we replace ``linearly independent'' by ``orthonormal''.}
\end{remark}

\section{A relationship between $\rho^{\mathbb H}({\mathbb H}^n)$ and $\rho^{\mathbb C}({\mathbb C}^{2n})$}

For a positive integer $m$, let $c_m^{\mathbb F}$ be the
${\mathbb F}$-James number mentioned in \cite{Jam58}. Then  $\rho^{\mathbb F}({\mathbb F}^n)$ is the largest $m$ such that $n$ is a multiple of $c_m^{\mathbb F}$.
$c_m^{\mathbb R}$ is computed by Adams in \cite{Ada62} and
is given by $c_m^{\mathbb R}= 2^{f(m)}$ where $f(m)$ is the number of
integers $q$ such that $0<q\leqslant m$ and $q\equiv0,1,2\mbox{ or } 4 \mbox{ mod } 8$.
In \cite{Ada62}, Adams also showed that $\rho^{\mathbb R}({\mathbb R}^n) = 8d+2^c-1$ where
$n = (2a+1)2^{4d+c}$ for some non-negative integers $a, d, c$ with $0\leqslant c \leqslant 3$.

Now, we show how to compute  $\rho^{\mathbb C}({\mathbb C}^n),
\rho^{\mathbb H}({\mathbb H}^n)$. For a given prime $p$, let $\nu_p(m)$ be the exponent of $p$ in the prime decomposition of $m$.
$c_m^{\mathbb C}$ is computed by Adams-Walker \cite{AW65}
and is given by $$\nu_p(c_m^{\mathbb C})=\max\{s+\nu_p(s):0\leqslant s \leqslant \lfloor \frac{m}{p-1}\rfloor\}.$$
$c_m^{\mathbb H}$ is computed by Sigrist-Suter \cite{SS73}
and is given by $\nu_2(c_m^{\mathbb H})=\max \{2m+1, 2s+\nu_2(s):0\leqslant s \leqslant m\},$ and for an odd prime $p$,
 $$\nu_p(c_m^{\mathbb H})=\max \{s+\nu_p(s):0\leqslant s \leqslant \lfloor \frac{2m}{p-1}\rfloor\}.$$
For $p\leqslant m+1$, let $$S_{m, p}=\{s:\lfloor \frac{m}{p-1}
\rfloor+\nu_p(\lfloor \frac{m}{p-1}\rfloor )-\lfloor \log_p(\frac{m}{p-1})\rfloor \leqslant s \leqslant \lfloor \frac{m}{p-1}\rfloor \},$$
and for $p > m+1$,  let $S_{m, p}=\{0\}$.
Notice that if $p\leqslant m+1$, then $|S_{m, p}|= \lfloor \log_p(m/(p-1))\rfloor - \nu_p(\lfloor m/(p-1)\rfloor ) + 1$. In the following lemma, we give refined formulas for computing
$\nu_p(c_m^{\mathbb C})$  and $\nu_p(c_m^{\mathbb H})$.

\begin{lemma} \label{lem:JamesNumbers}
Let $m$ be  a positive integer and $p$ be a prime number. Then
\begin{enumerate}
\item[\textup{(i)}] $\nu_p(c_m^{\mathbb C})=\max\{s + \nu_p(s):s\in S_{m,p}\}$.
\item[\textup{(ii)}] $\nu_2(c_m^{\mathbb H})=\max\{2m+1, 2s + \nu_2(s):s\in S_{m,2}\}$, and for an odd prime $p$, $\nu_p(c_m^{\mathbb H})=\max\{s + \nu_p(s):s\in S_{2m,p}\}$.
\end{enumerate}
\end{lemma}
\textbf{Proof.}
\textup{(i)} Let $\nu_p(c_m^{\mathbb C})=u + \nu_p(u)$ for some
$0\leqslant u \leqslant \lfloor m/(p-1)\rfloor$. Then
$\nu_p(u)\leqslant \lfloor \log_pu\rfloor \leqslant \lfloor\log_p(m/(p-1))\rfloor$. So,
$\lfloor m/(p-1)\rfloor+\nu_p(\lfloor m/(p-1)\rfloor )\leqslant u + \nu_p(u)\leqslant u + \lfloor\log_p(m/(p-1))\rfloor$.  Hence, $u\geqslant \lfloor m/(p-1)\rfloor+\nu_p(\lfloor m/(p-1)\rfloor - \lfloor\log_p(m/(p-1))\rfloor$.The result follows.
(ii) is similar to (i).\\

Let $p_i$ be the $i$th prime number, i.e.,
$p_1 = 2,p_2 = 3, p_3 = 5, \mbox{etc}$. Let $n = p_1^{t_1}\times \cdots \times p_r^{t_r}\times l$, where $t_i\geqslant 1$ for each $i=1,\ldots, r$ and $\nu_{p_{r+1}}(l) = 0$. For each $i = 1, \ldots, r$, let
$K_{n,i}^{\mathbb C}$ be the largest  $m\in\{0,\ldots,p_{r+1}-2\}$ such that
$t_i\geqslant \nu_{p_i}(c_m^{\mathbb C})$, and let $K_{n,i}^{\mathbb H}$
be the largest $m\in\{0,\ldots,(p_{r+1}-3)/2\}$ such that
$t_i\geqslant \nu_{p_i}(c_m^{\mathbb H})$. In the following theorem, we give direct formulas for computing $\rho^{\mathbb C}({\mathbb C}^n)$ and
$\rho^{\mathbb H}({\mathbb H}^n)$.

\begin{theorem} \label{theorem:FormulaComplexQuaternionicVectorFields}
Let $n = p_1^{t_1}\times \cdots \times p_r^{t_r}\times l$, where $t_i\geqslant 1$ for each $i=1,\ldots, r$ and $\nu_{p_{r+1}}(l) = 0$. Then
\begin{enumerate}
\item[\textup{(i)}]$\rho^{\mathbb C}({\mathbb C}^n) =
\min\{K_{n,i}^{\mathbb C}:i=1,\ldots,r\}$.
\item[\textup{(ii)}] $\rho^{\mathbb H}({\mathbb H}^n) =
\min\{K_{n,i}^{\mathbb H}:i=1,\ldots,r\}$.
\end{enumerate}
\end{theorem}
\textbf{Proof. }\textup{(i)} First, since $\nu_{p_{r+1}}(c_{p_{r+1}-1}^{\mathbb C}) = 1$ and $c_i^{\mathbb C}\mid c_{i+1}^{\mathbb C}$
for each $i\geqslant 1$ then $\rho^{\mathbb C}({\mathbb C}^n)\leqslant K_{n,i}^{\mathbb C}$ for each $i=1,\ldots,r$. On the other hand, if
$s = \min\{K_{n,i}^{\mathbb C}:i=1,\ldots,r\}$ then
$t_i\geqslant \nu_{p_i}(c_s^{\mathbb C})$ for each $i=1,\ldots,r$, and  hence $\rho^{\mathbb C}({\mathbb C}^n)\geqslant s$. The result follows. (ii) is similar to (i).

From Corollary~\ref{corollary:RelationshipNumberOfVectorFields}, we know that the number of linearly independent complex vector fields on $S({\mathbb C}^n)$ is at most half the  number of linearly independent real vector fields on $S({\mathbb R}^{2n})$, and the number  of linearly independent quaternionic vector fields on $S({\mathbb H}^{n})$ is at most half the number of linearly independent complex vector fields on $S({\mathbb C}^{2n})$. The main result of this
section is the following theorem , which gives an explicit relationship between
the number of linearly independent complex vector fields on $S({\mathbb C}^{2n})$
and the number  of linearly independent quaternionic vector fields on $S({\mathbb H}^{n})$.

\begin{theorem} \label{theorem:RelationshipNumberOfComplexQuaternionicVectorFields}
  $\rho^{\mathbb C}({\mathbb C}^{2n}) = 2\rho^{\mathbb H}({\mathbb H}^n)+d$
  where $d = 1\mbox{ or } 3$.
\end{theorem}
\textbf{Proof. } Let $\rho^{\mathbb H}({\mathbb H}^n)=m$. Then $m$ is the smallest integer such that $n$ is not a multiple of $c_{m+1}^{\mathbb H}$. From \cite{SS73}, $c_{m+1}^{\mathbb H}$ is either equal to $c_{2m+3}^{\mathbb C}/2$ or to $c_{2m+3}^{\mathbb C}$.
If $c_{m+1}^{\mathbb H}=c_{2m+3}^{\mathbb C}/2$, then  $2n$ is not a multiple of
$2c_{m+1}^{\mathbb H}=c_{2m+3}^{\mathbb C}$.
On the other hand, if $c_{m+1}^{\mathbb H}=c_{2m+3}^{\mathbb C}$ then $\nu_2(c_{m+1}^{\mathbb H})= \nu_2(c_{2m+3}^{\mathbb C})= 2m+3$. Hence, $2n$ is not a multiple of $c_{m+2}^{\mathbb H}$, because $c_{m+2}^{\mathbb H}\geqslant c_{m+1}^{\mathbb H}$, $\nu_p(2n)=\nu_p(n)$ for each $p\neq 2$,  and  $\nu_2(c_{m+2}^{\mathbb H})\geqslant 2m+5$, see Lemma~\ref{lem:JamesNumbers}. But, $c_{m+2}^{\mathbb H}$ is either equal to $c_{2m+5}^{\mathbb C}/2$ or to $c_{2m+5}^{\mathbb C}$. Thus, in all cases, $2n$ is not a multiple of $c_{2m+5}^{\mathbb C}$, and hence
$\rho^{\mathbb C}({\mathbb C}^{2n})< 2m+5$. Now, since  $c_{2k+1}^{\mathbb C}=c_{2k}^{\mathbb C}$ for each $k\geqslant 1$, see \cite{AW65}, then $\rho^{\mathbb C}({\mathbb C}^{2n})$ is odd. Consequently, by Corollary~\ref{corollary:RelationshipNumberOfVectorFields}, $\rho^{\mathbb C}({\mathbb C}^{2n})$
is either equal to  $2m+1$  or  to $2m+3$.\\

In the following table, we give the values of $\rho^{\mathbb R}({\mathbb R}^{4n})$, $\rho^{\mathbb C}({\mathbb C}^{2n})$, and $\rho^{\mathbb H}({\mathbb H}^{n})$ for some values of $n$.

\begin{table}[h]
\centering
\begin{tabular}{|c|c|c|c|c|c|c|c|}
\hline
$n$&
${1}$&
${2}$&
${4}$&
${6}$&
${12}$&
${24}$&
$1440$ \\
\hline
$\rho^{\mathbb R}({\mathbb R}^{4n})$&
3&
7&
8&
7&
8&
9&
15 \\
\hline
$\rho^{\mathbb C}({\mathbb C}^{2n})$&
1&
1&
1&
1&
3&
3&
5 \\
\hline
$\rho^{\mathbb H}({\mathbb H}^{n})$&
0&
0&
0&
0&
0&
1&
2 \\
\hline
\end{tabular}
\end{table}

\section{Construction of vector fields on spheres}
The construction of $\rho^{\mathbb R}({\mathbb R}^{n})$ linearly independent real vector
fields on $S({\mathbb R}^{n})$ is well understood, for constructions using
Clifford algebras see \cite{Zve68}, and for
constructions using combinatorial methods see the recent work of Ognikyan
\cite{Ogn08}. The situation is completely different with
complex and quaternionic vector fields; there is no explicit constructions
that gives two or more linearly independent complex vector fields on $S({\mathbb C}^{n})$ and there is no known construction that gives even a single
quaternionic vector field on $S({\mathbb H}^{n})$. In fact, the only known
complex vector field on $S({\mathbb C}^{n})$ is the one given
in Example~\ref{exa:ComplexVectorField}.

From Theorem~\ref{theorem:LinearlyIndependentVectorFields}, we know that one can use $m$ linearly independent complex (respectively, quaternionic) vector fields to obtain
$2m$ linearly independent real (respectively, complex) vector
fields. So, it is natural to ask if, in some way, it is possible to use
linearly independent real (respectively, complex) vector fields to build
some linearly independent complex or quaternionic (respectively, quaternionic) vector
fields.

In the following theorem, we give necessary and sufficient conditions on
linearly independent real (respectively, complex) vector fields to be
linearly independent complex or quaternionic (respectively, quaternionic)
vector fields. For simplicity, let $r_{{\mathbb C},i}=r_{\mathbb C}\circ \alpha_i\circ
 r_{\mathbb C}^{-1}$, $c_{{\mathbb H},j}=r_{\mathbb C}\circ \alpha_j\circ
 r_{\mathbb H}^{-1}$, and for each $t\in \{i,j,k\}$,  let $r_{{\mathbb H},t}=r_{\mathbb H}\circ \alpha_t\circ r_{\mathbb H}^{-1}$.

\begin{theorem} \label{theorem:ConstructionOfVectorFields}
\begin{enumerate}
\item[\textup{(i)}] Suppose $u_1,\ldots,u_m$ are linearly independent real vector
fields on $S({\mathbb R}^{2n})$. Then
$r_{\mathbb C}^{-1}\circ u_1 \circ r_{\mathbb C},\ldots,
r_{\mathbb C}^{-1}\circ u_m \circ r_{\mathbb C}$ are linearly independent  complex vector field on $S({\mathbb C}^n)$ if and only if
$u_1,\ldots,u_m,r_{{\mathbb C},i}\circ u_1,\ldots,r_{{\mathbb C},i}\circ u_m$ are linearly independent real vector fields on $S({\mathbb R}^{2n})$.
\item[\textup{(ii)}]  Suppose $v_1,\ldots,v_m$ are linearly independent complex vector
fields on $S({\mathbb C}^{2n})$. Then
$c_{\mathbb H}^{-1}\circ v_1 \circ c_{\mathbb H},\ldots,
c_{\mathbb H}^{-1}\circ v_m \circ c_{\mathbb H}$ are linearly independent quaternionic vector field on $S({\mathbb H}^n)$ if and only if
$v_1,\ldots,v_m,c_{{\mathbb H},j}\circ v_1,\ldots,c_{{\mathbb H},j}\circ v_m$ are linearly independent complex vector fields on $S({\mathbb C}^{2n})$.
\item[\textup{(iii)}]Suppose $u_1,\ldots,u_m$ are linearly independent real vector
fields on $S({\mathbb R}^{4n})$. Then
$r_{\mathbb H}^{-1}\circ u_1 \circ r_{\mathbb H},\ldots,
r_{\mathbb H}^{-1}\circ u_m \circ r_{\mathbb H}$ are linearly independent  quaternionic vector field on $S({\mathbb H}^n)$ if and only if
$u_1,\ldots,u_m,r_{{\mathbb H},t}\circ u_1,\ldots,r_{{\mathbb H},t}\circ u_m$, where
$t\in \{i,j,k\}$,  are linearly independent real vector fields on
$S({\mathbb R}^{4n})$.
\end{enumerate}
\end{theorem}
\textbf{Proof.} Follows directly from Theorem~\ref{theorem:LinearlyIndependentVectorFields}.\\

In \cite{Bec72}, Becker solved several important cases of
the equivariant real vector fields problem, both the maximal number and the
construction, on spheres with free group action by using the known Clifford
algebras constructions of the non-equivariant real vector fields on spheres.
As noted in \cite{Obi06} and \cite{Ond01}, both sides of the equivariant complex and quaternionic vector fields problem on spheres with free group action is still completely
open. By using a method similar to that used by Becker, the construction of
the non-equivariant complex and quaternionic vector fields on spheres might
lead to the solution of the equivariant complex and quaternionic vector
fields problem on spheres with free group action.

\end{document}